\documentclass[a4paper,12pt]{article}
\usepackage{latexsym}
\usepackage{amsfonts}
\usepackage{xypic}
\newtheorem{theorem}{Theorem}[section]
\newtheorem{corollary}[theorem]{Corollary}
\newtheorem{definition}{Definition}
\newtheorem{example}[theorem]{Example}
\newtheorem{remark}[theorem]{Remark}

\newtheorem{proposition}[theorem]{Proposition}

\newtheorem{observation}[theorem]{Observation}
\title{Idealization of some weak separation axioms\thanks{1991
Math.\ Subject Classification --- Primary: 54A05, 54D10;
Secondary: 54D30, 54H05. \protect\newline Key words and phrases ---
topological ideal, $T_0$, \two, \three, \four.}}
\author{Francisco G. Arenas\footnote{Research supported
by DGICYT grant PB95-0737.}, Julian Dontchev\footnote{Research
supported partially by the Ella and Georg Ehrnrooth Foundation
at Merita Bank, Finland.} and Maria Luz Puertas}
\begin{document}
\baselineskip=20pt plus 1pt minus 1pt
\newcommand{\xti}{$(X,\tau,{\cal I})$}
\newcommand{\ysj}{$(Y,\sigma,{\cal J})$}
\newcommand{\ti}{$T_{\cal I}$}
\newcommand{\two}{$T_{\frac{1}{2}}$}
\newcommand{\three}{$T_{\frac{1}{3}}$}
\newcommand{\four}{$T_{\frac{1}{4}}$}
\newcommand{\fxy}{$f \colon (X,\tau) \rightarrow (Y,\sigma)$}
\newcommand{\fxyi}{$f \colon (X,\tau,{\cal I}) \rightarrow
(Y,\sigma,{\cal J})$}
\maketitle
\begin{abstract}
An ideal is a nonempty collection of subsets closed under
heredity and finite additivity. The aim of this paper is to unify
some weak separation properties via topological ideals. We
concentrate our attention on the separation axioms between $T_0$
and \two. We prove that if \xti\ is a semi-Alexandroff $T_{\cal
I}$-space and $\cal I$ is a $\tau$-boundary, then $\cal I$ is
completely codense.
\end{abstract}

\section{Introduction}\label{s1}

Regardless of the fact that for many years ideals had their
significant impact in research in Topology it was probably the
four articles of Hamlett and Jankovi\'{c} \cite{HJ1,HJ2,JH1,JH2},
which appeared almost a decade ago, that initiated the
application of topological ideals in the generalization of most
fundamental properties in General Topology.

Throughout the 90's several topological properties such as
covering properties, connectedness, resolvability, extremal
disconnectedness and submaximality have been generalized via
topological ideals. Contributions in the field are due to (in
alphabetical order) Abd El-Monsef, Ergun, Ganster, Hamlett,
Jankovi\'{c}, Lashien, Maki, Nasef, Noiri, Rose and Umehara.

Probably only separation axioms have been neglected from `ideal
point of view'. In this paper, we will consider only weak
separation axioms -- the ones between $T_0$ and \two. Recall that
a topological space $(X,\tau)$ is called a {\em \two-space} if
every singleton is open or closed. The importance of the
separation axiom \two\ is probably given by the following
(perhaps well-known) result: Every minimal \two-space is
compact and connected.

In Digital Topology \cite{KR1} several spaces that fail to be
$T_1$ are important in the study of the geometric and topological
properties of digital images \cite{KR1,KK1}. Such is the case
with the major building block of the digital n-space -- the {\em
digital line} or the so called {\em Khalimsky line}. This is the
set of the integers, $\mathbb Z$, equipped with the
topology $\cal K$, generated by ${\cal G}_{\cal K} = \{ \{ 2n-1,
2n, 2n+1 \} \colon n \in {\mathbb Z} \}$. Although the digital
line is not a $T_1$-space, it satisfies the separation axiom
\two. This probably indicates that further knowledge (from more
general point of view) of the behavior of topological spaces
satisfying separation axioms below $T_1$ is required, which
indeed is the intention of the present paper.

A {\em topological ideal} is a nonempty collection of subsets of
a topological space $(X,\tau)$, which is closed under heredity
and finite additivity. Except the trivial ideals, i.e. the
minimal ideal $\{ \emptyset \}$ and the maximal ideal $\cal P$,
the following collections of sets form important ideals on any
topological space $(X,\tau)$: the finite sets $\cal F$, the
countable sets $\cal C$, the closed and discrete sets $\cal CD$,
the nowhere dense sets $\cal N$, the meager sets $\cal M$, the
scattered sets $\cal S$ (only when $X$ is $T_0$), the bounded
sets $\cal B$, the relatively compact sets $\cal R$, the
hereditarily compact (resp.\ Lindel\"of) sets $\cal HK$ (resp.\
$\cal HL$), the $S$-bounded sets ${\cal SB}$ \cite{K1} and (in
the real line) the Lebesgue null sets $\cal L$.

An {\em ideal topological space} \xti\ is a topological space
$(X,\tau)$ and an ideal $\cal I$ on $(X,\tau)$. For an ideal
topological space \xti\ and a subset $A \subseteq X$,
$A^{*}({\cal I}) = \{x \in X: U \cap A \not\in {\cal I}$ for
every $U \in {\tau} (x)\}$ is called the {\em local function\/}
of $A$ with respect to $\cal I$ and $\tau$ \cite{K2} (where
$\tau(x)$ is the open neighborhood filter at $x$). We simply
write $A^{*}$ instead of $A^{*}({\cal I})$ in case there is no
chance for confusion.

Note that ${\rm Cl}^{*}(A) = A \cup A^{*}$ defines a Kuratowski
closure operator for a topology $\tau^{*}({\cal I})$ (also
denoted by $\tau^{*}$ when there is no chance for confusion),
finer than $\tau$. A basis $\beta({\cal I},\tau)$ for
$\tau^{*}({\cal I})$ can be described as follows: $\beta({\cal
I},\tau) = \{ U \setminus I \colon U \in \tau$ and $I \in {\cal
I}\}$. In general, $\beta$ is not always a topology \cite{JH1}.

\section{Weak separation axioms below $T_1$ via ideals}\label{s2}

\begin{definition}\label{d1}
{\em An ideal topological space \xti\ is called a {\em $T_{\cal
I}$-space} if for every subset $I \in {\cal I}$ of $X$ and every
$x \not\in I$, there exists a set $A_x$ containing $x$ and
disjoint from $I$ such that $A_x$ is open or closed.}
\end{definition}

\begin{proposition}\label{p4}
If a topological space $X$ is $T_{\mathcal{J}}$, then it is a
$T_{\mathcal{I}}$-space for each ideal $\mathcal{I}\subseteq
\mathcal{J}$.
\end{proposition}

{\em Proof.} Let $X$ be a $T_{\mathcal{J}}$-space. Suppose that
$\mathcal{I}$ is an ideal with $\mathcal{I}\subseteq{J}$ and let
$I\in \mathcal{I}$ and $x\notin I$. Then $I \in \mathcal{J}$ and
there exists a set $A_x$ which is open or closed such that $x\in
A_x$ and $A_x\cap I=\emptyset$, so $X$ is
$T_{\mathcal{I}}$-space. $\Box$

\begin{observation}\label{o1}
{\rm (i) Every \two-space is a $T_{\cal I}$-space.

(ii) Every space is a $T_{\cal CD}$-space.

(iii) If $\mathcal{I}$ is an ideal on a set $X$, the topological
space $X$ with the family of closed sets $\mathcal{I}\cup\{ X\}$,
is a $T_{\mathcal{I}}$-space. Note that using the notation from
the preceding section, this topology on $X$ can be described as
$\tau_{T}^*({\cal I})$, where $\tau_{T}$ is the trivial topology
in $X$, and in this case the basis $(\beta{\cal I},\tau_{T})$ is
in fact the whole topology. On the other hand, if this
topological space is a $T_{\mathcal{J}}$-space for another ideal
${\cal J}$, then the following relation between both ideals
${\cal I}$ and ${\cal J}$ are obtained: if $F\in {\cal
J}\setminus {\cal I}$ then $F$ is the intersection of some open
sets of $X$. To see this, suppose that $F\in {\cal J}\setminus
{\cal I}$, then for each $x\notin F$ there exists a closed set
$A_x$ disjoint from $F$ that contains $x$ (note that if
$A_x = X\setminus I$, with $I\in {\cal I}$, is a non-empty open
set disjoint from $F$, then $F\subseteq I$ and $F\in {\cal I})$.
So, $F=\cap_{x\notin F}(X \setminus A_x)$ where $X \setminus A_x$
is open for each $x \notin F$.}
\end{observation}

In \cite{Nc1}, Newcomb defined an ideal $\cal I$ on a space
$(X,\tau)$ to be a {\em $\tau$-boundary\/} if $\tau \cap {\cal
I} = \{ \emptyset \}$. A topological space $(X,\tau)$ is called
{\em semi-Alexandroff} \cite{ADG2} if any intersection of open
sets is semi-open, where a {\em semi-open} set is a set which can
be placed between an open set and its closure. Complements of
semi-open sets are called {\em semi-closed}.

In connection with our next result, we recall that a set $A$ is
called {\em locally dense} \cite{CM1} or {\em preopen} if $A
\subseteq {\rm Int} \overline{A}$. The collection of all preopen
subsets of a topological space $(X,\tau)$ is denoted by $PO(X)$.
An ideal $\cal I$ on a space \xti\ is called {\em completely
codense} \cite{DGR1} if $PO(X) \cap {\cal I} = \{ \emptyset \}$.
Note that if $({\mathbb R},\tau)$ is the real line with the usual
topology, then $\cal C$ is codense but not completely codense.
It is shown in \cite{DGR1} that an ideal $\cal I$ is completely
codense on $(X,\tau)$ if and only if ${\cal I} \subseteq {\cal
N}$, i.e.\ if each member of $\cal I$ is nowhere dense.

\begin{proposition}\label{al1}
If \xti\ is a semi-Alexandroff $T_{\cal I}$-space and $\cal I$
is a $\tau$-boundary, then $\cal I$ is completely codense.
\end{proposition}

{\em Proof.} Let $I \in {\cal I}$. If $I = X$, i.e., if ${\cal
I} = {\cal P}$, we have $X = \emptyset$, since $\cal I$ is a
$\tau$-boundary and we are done. Assume next that $I$ is a proper
subset of $X$. Since \xti\ is a $T_{\cal I}$-space, for every
point $x \not \in I$, we can find a set $A_x$ such that $x \in
A_x$, $A_x \cap I = \emptyset$ and $A_x$ is open or closed. Set
$U = \cup \{ A_x \colon A_x$ is open$\}$ and $V = \cup \{ A_x
\colon A_x$ is not open$\}$. Clearly, $U \in \tau$ and $V$ is
semi-closed, since $(X,\tau)$ is a semi-Alexandroff space. Thus,
$X \setminus I$ is the union of an open and a semi-closed set.
Hence, $I$ is the intersection of a semi-open and a semi-closed
set or equivalently the union of an open set $W$ and a nowhere
dense set $N$ (such sets are called simply-open). Note that $W
= \emptyset$, since $\cal I$ is a $\tau$-boundary and $W \in
{\cal I}$ due to the heredity of the ideal $\cal I$. This shows
that $I = N$, i.e., ${\cal I} \subseteq {\cal N}$. $\Box$

Recall that a nonempty topological space $(X,\tau)$ is called
{\em resolvable} \cite{H1} (resp.\ {\em $\cal I$-resolvable}
\cite{DGR1}) if $X$ is the disjoint union of two dense (resp.\
$\cal I$-dense) subsets, where a subset $A$ of a topological
space \xti\ is called {\em $\cal I$-dense} \cite{DGR1} if every
point of $X$ is in the local function of $A$ with respect to
$\cal I$ and $\tau$, i.e.\ if $A^{*} ({\cal I}) = X$.

\begin{corollary}
If \xti\ is a semi-Alexandroff $T_{\cal I}$-space such that $\cal
I$ is a $\tau$-boundary, then resolvability of $(X,\tau)$ implies
automatically the $\cal I$-resolvability of \xti.
\end{corollary}

Several already defined separation axioms are either implied or
equivalent to the `ideal separation axiom' in Definition~\ref{d1}
when a certain ideal is considered.

Recall that a topological space $(X,\tau)$ is called a {\em
$T_{{\frac{1}{4}}}$-space} \cite{ADG1} (resp.\ a {\em
$T_{{\frac{1}{3}}}$-space} \cite{ADP1}) if for every finite
(resp.\ compact) subset $S$ of $X$ and every $x \not\in S$
there exists a set $S_x$ containing $F$ and disjoint from $\{ x
\}$ such that $A_x$ is open or closed.

\begin{proposition}\label{p1}
Let $(X,\tau)$ be a topological space. Then the following
conditions are valid:

{\rm (i)} $(X,\tau)$ is a \four-space if and only if \xti\ is a
$T_{\cal F}$-space.

{\rm (ii)} Every \three-space is a $T_{\cal HK}$-space.

{\rm (iii)} $(X,\tau)$ is a \two-space if and only if \xti\ is
a $T_{\cal P}$-space.

{\rm (iv)} Every $T_{{\cal I}(K)}$-space is a \three-space, where
${\cal I}(K)$ is the ideal formed by the subsets of the compact
sets of the topological space.
\end{proposition}

{\em Proof.} (i) and (ii) follow directly from the definitions.

(iii) Assume that $(X,\tau)$ is a two-space, let $P$ be an
arbitrary subset of $X$ and let $x \not\in P$. Then $\{ x \}$ is
open or closed, which points out that \xti\ is a
$T_{\cal P}$-space. Conversely, if \xti\ is a $T_{\cal P}$-space
and $x \in X$, then there exists a set $A_x$, which is open or
closed such that $x \in A_x$ and $A_x \cap X \setminus \{ x \}
= \emptyset$. Clearly, $\{ x \} = A_x$. This shows that every
singleton of $(X,\tau)$ is open or closed or equivalently
$(X,\tau)$ is a \two-space.

{\rm (iv)} Let $X$ be a $T_{{\mathcal I}(K)}$-space, $K$ a
compact subset of $X$ and $x\notin K$. Then $K\in {\mathcal
I}(K)$, so there exists a subset $A_x$ of $X$ which is either
open or closed such that $x\in A_x$ and $A_x\cap K=\emptyset$,
so $X$ is a \three-space. $\Box$

\begin{example}
{\em We give an example of a $T_{\cal HK}$-space which is not a
\three-space. We topologize the real line $\mathbb R$ by
declaring the following non-trivial open sets: 

(1) all singletons except $\{1\}$ and $\{2\}$; 

(2) all cofinite sets containing $2$ but not $1$;

(3) all cofinite sets containing the set $\{ 1,2 \}$.

We show that this is a $T_{\cal HK}$-space. Let $A \subseteq
{\mathbb R}$ be hereditarily compact and let $x \not\in A$. If
$x \not= 2$, then $\{ x \}$ is open or closed and we are done.
If $x = 2$, then $1 \not\in A$, in which case $A$ is open, or $1
\in A$, in which case $A \setminus \{ 1 \}$ is finite for $A$ is
hereditarily compact and hence $A$ is closed. 

Next we observe that the space is not a \three-space. For note
that $A \setminus \{ 2 \}$ is compact but $\{ 2 \}$ is neither
open nor closed.}
\end{example}

If $X$ is a set, a {\em subideal} (or an {\em h-family}) $\cal
I$ on $X$, is a family of subsets of $X$ which is closed under
the subset relation. A family which is closed under finite
additivity is called an {\em FA-family}. Sometimes, in the
definition of $T_{\cal I}$-spaces, we might not require that the
family $\cal I$ is an ideal, we take the more general approach
of considering subideals and FA-families. The FA-family formed
by Lindel\"of subsets of a given topological space will be
denoted by ${\cal I}(L)$.

A set $A \subseteq (X,\tau)$ is called {\em discretely finite}
(= df-set) (resp.\ {\em discretely countable} (=dc-set))
\cite{DGR1} if for every point $x \in A$, there exists $U \in
\tau$ containing $x$ such that $U \cap A$ is finite (resp.\
countable). The subideals of all df-sets (resp.\ dc-sets) will
be denoted by ${\cal DF}$ (resp.\ ${\cal DC}$).

Using Proposition~\ref{p4}, which is valid for pairs of subideals
and pairs of FA-families, we obtain the following diagram where
the relations between different classes of $T_{\cal I}$-spaces
are shown.

$$
\diagram
\text{$T_1$} \rto \dto & \text{$T_{\cal HL}$} \rto \dto &
\text{$T_{{\cal I}(L)}$} \dto \rto & \text{$T_{\cal C}$} \dto
\rto & \text{$T_{\cal DC}$} \dto \rto & \text{$T_0$} \dto
\\ \text{$T_{\cal P}$=$T_{\frac{1}{2}}$} \rto & \text{$T_{\cal
HK}$} \rto & \text{$T_{{\cal I}(K)}$} \rto & \text{$T_{\cal F}$}
\rto & \text{$T_{\cal DF}$} \urto \rto & \text{$T_{\emptyset}$}
\enddiagram
$$

\begin{example}
{\em A \four-space which is not a $T_{\mathcal{HK}}$-space. The
topological space in \cite{ADP1}, Example 3.1 is a \four-space.
This is the set $X$ of non-negative integers with the topology
whose open sets are those which contain $0$ and have finite
complement (so closed sets are the finite sets that do not
contain $0$). Each subset of $X$ is compact, then ${\mathcal
HK}=\mathcal{P}$ and $X$ is not a $T_{\mathcal{HK}}$-space
because it is not a \two-space.}
\end{example}

Recall that an {\em $\alpha$-space} (or a {\em nodec space}) is
a space where $\tau^{*}({\cal N}) = \tau$. Note that
$\alpha$-spaces need not be even $T_0$.

\begin{remark}
{\rm 1) Every $\alpha$-space $(X,\tau)$ is a $T_{\cal N}$-space:
Let $N \in {\cal N}$ and $x \not\in N$. Since $X$ is an
$\alpha$-space, then every nowhere dense subset is closed (and
discrete). Thus $N$ is closed, which shows that $X$ is a $T_{\cal
N}$-space.

(2) It is not difficult to find an example of a $T_{\cal
N}$-space which is not an $\alpha$-space. Note that even a metric
space (for example the real line with usual topology) need not
be an $\alpha$-space.}
\end{remark}

A subset $S$ of a space $(X,{\tau},{\cal I})$ is a topological
space with an ideal ${\cal I}_{S} = \{ I \in {\cal I} \colon I
\subseteq S \} = \{ I \cap S \colon I \in {\cal I} \}$ on $S$
\cite{JD1}.

\begin{proposition}\label{p3}
Every subspace of a \ti-space is a $T_{{\cal I}_A}$-space.
\end{proposition}

{\em Proof.} Let \xti\ be a \ti-space and let $A \subseteq X$.
If $J \in {\cal I}_{A}$ is a subset of $(A,\tau|A,{\cal I}_{A})$
and $x \in A \setminus J$, then there exists $I \in {\cal I}$
such that $I \cap A = J$. Since \xti\ is a \ti-space, there
exists a set $B \supseteq I$ such that either $B \in \tau$ or $X
\setminus B \in \tau$ and $x \not \in B$. Clearly, $B \cap A$ is
open or closed in the subspace $(X,\tau|A)$. This shows that
$(A,\tau|A,{\cal I}_{A})$ is a $T_{{\cal I}_A}$-space. $\Box$

\begin{proposition}\label{mp3}
Every topology finer than a $T_{\mathcal{HK}}$-topology is a
$T_{\mathcal{HK}}$-topology.
\end{proposition}

{\em Proof.} Let $\tau^* \subseteq \tau$ be two topologies on a
set $X$. Then, $\mathcal{HK}_{\tau}\subseteq \mathcal{HK}_{\tau
^*}$: let $H \in \mathcal{HK}_{\tau}$ be and $F \subseteq H$,
then $F$ is $\tau$-compact, so it is $\tau^*$-compact and $H \in
\mathcal{HK}_{\tau^*}$.

Suppose that $(X,\tau^*)$ is a $T_{\mathcal{HK}_{\tau^*}}$-space.
Then $(X,\tau)$ is a $T_{\mathcal{HK}_{\tau^*}}$-space: let $H
\in \mathcal{HK}_{\tau^*}$ and $x \notin H$, then there exists
$A_x$ which is $\tau^*$-closed or $\tau ^*$-open such that $x
\notin A_x$ and $A_x \cap H = \emptyset$. Then $A_x$ is
$\tau$-closed or $\tau$-open and so $(X,\tau)$ is a
$T_{\mathcal{HK}_{\tau^*}}$-space. Now, using that
$\mathcal{HK}_{\tau} \subseteq \mathcal{HK}_{\tau^*}$ and
Proposition~\ref{p4}, we obtain that $(X,\tau)$ is a
$T_{\mathcal{HK}_{\tau}}$-space. $\Box$

\begin{remark}
{\rm Note that if ${\cal I}$ is an ideal of a topological space
$X$ that does not depend on the topology, for example
$\mathcal{F}$ or $\mathcal{C}$, then with a proof similar to the
one below every topology in $X$ finer that a $T_{\cal
I}$-topology is a $T_{\cal I}$-topology.}
\end{remark}

A property $\cal P$ of an ideal topological space is called {\em
I-topological} if \ysj\ has the property $\cal P$ whenever \xti\
has the property $\cal P$ and \fxyi\ is an ${\cal
IJ}$-homeomorphism, i.e., \fxy\ is a homeomorphism and $f({\cal
I}) = {\cal J}$. It is well-known that the property \two\ is
topological. The following proposition is an ideal generalization
of that result.

\begin{proposition}
The property \ti\ is I-topological.
\end{proposition}

{\em Proof.} Assume that \xti\ is a \ti-space and that \fxyi\ is
an ${\cal IJ}$-homeomorphism. Let $J \in {\cal J}$ and let $y \in
Y$ such that $y \not\in J$. Since \fxyi\ is an I-homeomorphism,
then there exists $x \in X$ and $I \in {\cal I}$ such that $f(I)
= J$, $y = f(x)$ and $x \not\in I$. Then, find a set $A_x
\subseteq X$ such that $x \in A_x$, $A_x \cap I = \emptyset$ and
$A_x$ is open or closed (such choice is possible for \xti\ is a
\ti-space). Since $f$ is open and closed, then $f(A_{x})$ is open
or closed subset of $(Y,\sigma)$. Moreover, $y \in f(A_{x})$ and
$f(A_{x}) \cap J = \emptyset$. This shows that \ysj\ is a
$T_{\cal J}$-space. $\Box$

{\bf Question.} How could we ideally extend the following result:
Every minimal \two-space is compact and connected.

\
\
\begin{center}
Area of Geometry and Topology\\Faculty of
Science\\Universidad de Almer\'{\i}a\\04071
Almer\'{\i}a\\Spain\\e-mail: {\tt farenas@ualm.es}
\end{center}
\
\begin{center}
Department of Mathematics\\University of Helsinki\\PL 4,
Yliopistonkatu 5\\00014 Helsinki\\Finland\\e-mail: {\tt
dontchev@cc.helsinki.fi}, {\tt
dontchev@e-math.ams.org}\\http://www.helsinki.fi/\~{}dontchev/
\end{center}
\
\begin{center}
Area of Geometry and Topology\\Faculty of
Science\\Universidad de Almer\'{\i}a\\04071
Almer\'{\i}a\\Spain\\e-mail: {\tt mpuertas@ualm.es}
\end{center}
\
\
\end{document}